\documentclass[12pt]{article}
\usepackage{amssymb,amsfonts,amsmath, amsthm}
\usepackage{adjustbox}
\usepackage{graphicx}
\usepackage{enumerate}
\usepackage{natbib}
\usepackage{url} 
\usepackage{tikz}

\newcommand{\blind}{1}

\newtheorem{cnj}{Conjecture}
\newtheorem{cor}{Corollary}

\newtheorem{rmk}{Remark}
\newtheorem{thm}{Theorem}

\begin{document}

\def\spacingset#1{\renewcommand{\baselinestretch}%
{#1}\small\normalsize} \spacingset{1}


\if1\blind
{
  \title{\bf More Power by Using Fewer Permutations}
  \author{Nick W. Koning\footnote{n.w.koning@ese.eur.nl. P.O. Box 1738, 3000 DR Rotterdam, The Netherlands.}\hspace{.2cm}\\
    Econometric Institute, Erasmus University Rotterdam}
  \maketitle
} \fi

\if0\blind
{
  \bigskip
  \bigskip
  \bigskip
  \begin{center}
    {\LARGE\bf More Power by using Fewer Permutations}
\end{center}
  \medskip
} \fi

\bigskip
\begin{abstract}
	It is conventionally believed that a permutation test should ideally use all permutations.
	If this is computationally unaffordable, it is believed one should use the largest affordable Monte Carlo sample or (algebraic) subgroup of permutations.
	We challenge this belief by showing we can sometimes obtain dramatically more power by using a tiny subgroup.
	As the subgroup is tiny, this simultaneously comes at a much lower computational cost.
	We exploit this to improve the popular permutation-based Westfall \& Young MaxT multiple testing method.
	We study the relative efficiency in a Gaussian location model, and find the largest gain in high dimensions.
	
\end{abstract}

\noindent%
{\it Keywords:}  group invariance test, maxT method, permutation test, randomization test, subgroup.
\vfill

\newpage

\section{Introduction}\label{sec:intro}
	Permutation- and more general ``group invariance''-based testing methods are fundamental tools in statistics.
	These methods date back to \citet{fisher1935} and have long been popular in causal inference and tests for independence, but also underly modern methodology such as large-scale permutation-based multiple testing methods \citep{westfall1993resampling, tusher2001significance} and conformal prediction \citep{shafer2008tutorial}.
	Moreover, even standard statistical tools such as the $t$-test can be interpreted as group invariance tests \citep{lehmann1949theory, koning2023more}.
	
	Fundamentally, these methods test the null hypothesis that the law of a random variable $X$ is invariant under a compact \emph{group} $\mathcal{G}$: 
	\begin{align}\label{eq:null}
		H_0 : X \overset{d}{=} G X, \text{ for all } G \in \mathcal{G}.
	\end{align}
	Equivalently, we can write $H_0 : X \overset{d}{=} \overline{G}X$, where  $\overline{G}$ is uniformly (Haar) distributed on $\mathcal{G}$ and independent of $X$, against the alternative that it is not invariant under $\mathcal{G}$.
	Instead of a direct interest in testing invariance, it is often the case that a null hypothesis implies invariance, so that rejecting this invariance also rejects the null hypothesis of interest.
	
	Observing a realization $x$ of $X$, the standard construction of a \emph{group invariance test} is as follows.
	One first chooses a test statistic $T$ that takes on extreme values under alternatives of interest.
	Then, the test rejects if its realization $T(x)$ exceeds a critical value that is an appropriate quantile of its \emph{reference distribution}: the distribution of $T(\overline{G}x)$ when $x$ is considered fixed.
	Such a test is exact up to discreteness of the reference distribution, and uniformly most powerful against simple alternatives if $T$ is chosen to be equal to the density of the alternative \citep{lehmann1949theory, koning2023online}.
	
\begin{figure}[!htb]
		\input{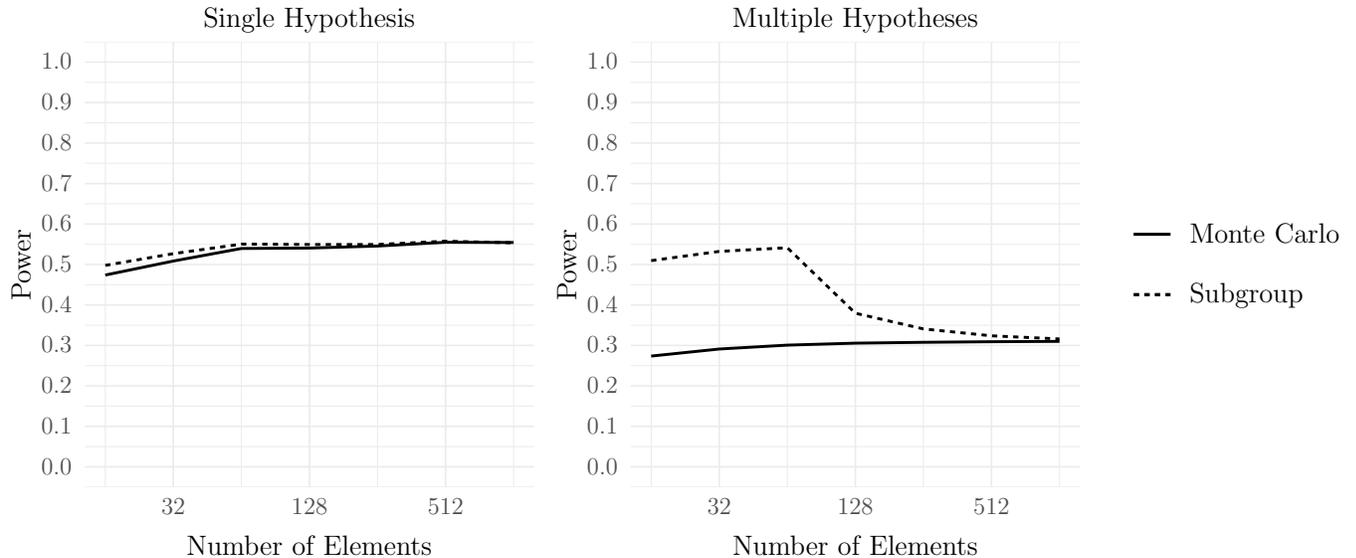}
		\caption{The power of subgroup (dotted) and Monte Carlo (solid) based methods, for single (left) and multiple (right) testing with the maxT method. The power is simulated based on 10\,000 repetitions for the one-sided testing of the location of one (left) and 1000 (right) locations in a Gaussian location model with $n = 32$ observations, all means equal to .3 (left) and .7 (right) and variance 1, at $\alpha = 1/16 = .0625$. The tests used are a standard sign-flipping test based on the sample mean as test statistic (left), and the maxT method based on 1000 of such sign-flipping tests (right). The Monte Carlo methods (solid) are based on the number of draws from the sign-flipping group indicated on the horizontal axis. The subgroup methods have a number of elements indicated on the horizontal axis, and are based on the best subgroups from the R-package NOSdata,  in such a way that the smaller subgroups are nested inside the larger ones.} 
		\label{fig:intro}
	\end{figure}

	Unfortunately, the group $\mathcal{G}$ of interest is often huge, making it computationally infeasible to use it in its entirety.
	As a result, it is near universal practice to approximate the reference distribution by a Monte Carlo approach that relies on drawing many times from the group \citep{eden1933validity, dwass1957, hemerik2018exact}.\footnote{In fact, this practice has seemingly become so universal that group invariance tests have grown to be nearly synonymous with randomization and are sometimes referred to as \emph{randomization tests}. See \citet{onghena2018randomization, hemerik2021another, zhang2023what, hemerik2023term} for a discussion on the relationship between randomization and group invariance tests, their history and terminology.}
	Such a Monte Carlo method comes with a trade-off: using more draws typically yields more power and higher replicability, but comes at a larger computational cost \citep{dwass1957, hope1968simplified}.
	This power-computation trade-off is visible in the left and right panel of Figure \ref{fig:intro}, where the solid line shows that the power of the single and multiple testing method displayed increases as the number of randomly drawn elements increases.

	\citet{chung1958randomization} proposed an alternative approach, in which they use the reference distribution based on a compact sub\emph{group} of $\mathcal{G}$.
	As invariance under a group implies invariance under its subgroups, such an approach still yields a valid test, but may be more computationally efficient if the subgroup is small.
	Recently, \citet{koning2023more} proposed to \emph{strategically select} the subgroup based on the choice of test statistic and alternative of interest.\footnote{We would like to stress that the selection of subgroups is also used in the conditional randomization literature. However, the goal there is to select a subgroup of permutations that guarantees (approximate) size control. In contrast, \citet{koning2023more} assume that a group that guarantees size control is available, and they subsequently select a subgroup that yields a particularly high powered test for its number of elements.}
	They find that this can more accurately approximate the reference distribution, leading to more power than the Monte Carlo approach for the same number of elements.
	Moreover, they empirically find that a larger subgroup yields more power, and recommend choosing the largest computationally affordable subgroup.
	An illustration of this can be seen in the left panel of Figure \ref{fig:intro}.
	There, the dotted line shows that the power of the subgroup-based method grows with the size of the subgroup, slightly faster than the power of the Monte Carlo-based method.
	
	\subsection{Contributions}	
		The main contribution of this paper is the observation that the power of group invariance based testing methods is not necessarily monotone in the number of elements of the subgroup, even if the subgroups are nested.
		This is illustrated in the right panel of Figure \ref{fig:intro}, where a method based on a subgroup with just 16 to 64 elements substantially outperforms the methods based on larger subgroups, as well as the Monte Carlo method based on a large number of draws.
		Contrary to popular belief, this implies one should not necessarily aspire to use the entire group.
		Indeed, if we would want to test invariance under the group with 1024 elements, then Figure \ref{fig:intro} suggests we would be better off using one of its subgroups with 64 elements, or even just 16 elements.
		The practical consequence of our finding is that we can obtain significantly \emph{more} power at a considerably \emph{lower} computational cost.
	
		The reason this finding was highly surprising to us, is that this seems to break a no-free-lunch principle: the subgroup-based method seemingly only exploits invariance under the subgroup, which is weaker than invariance under the entire group (and also weaker than invariance under larger subgroups that nest it).
		By this reasoning, a subgroup based method should not be expected to outperform the method based on the entire group.
		In fact, for a sequence of nested subgroups, we would expect the power to increase monotonically in the number of elements of the subgroup.
		
		We resolve this paradox as follows.
		While the use of an \emph{arbitrary} subgroup does not exploit invariance under the entire group, the \emph{selection} of the subgroup does exploit the fact that the null distribution is invariant under every subgroup of the entire group.
		This implies that a method based on a strategically selected subgroup \textit{does} exploit the invariance under the entire group.
		Moreover, by selecting the subgroup to have good power properties we can introduce knowledge about the alternative and test statistic, and thereby obtain higher power than a method based on the entire group.
		
		
		We apply our ideas in a high-dimensional location model for the group invariance-based \citet{westfall1993resampling} MaxT multiple testing method, which is perhaps the most popular multiple testing method for controlling the familywise error rate.
		We show the ideal subgroups rely on an orthogonality condition that can only hold if the subgroup has at most $n$ (for a two-sided hypothesis) or $2n$ (for a one-sided hypothesis) elements, which explains the drop in power in Figure \ref{fig:intro} as the subgroup grows beyond $2n = 64$ elements.
		The maxT method was also featured in one of the simulation studies of \citet{koning2023more} showing promising power results, but we find the subgroups they used were far too large.	
		
		Moreover, we analyze in which settings the power difference is largest, by studying the power of the subgroup- and entire group-based methods in a standard Gaussian location model, and comparing the two in terms of relative efficiency.
		Assuming all $p$ hypotheses are false, our results suggest that both tests are consistent for signals $\mu \succ n^{-1/2}\log^{1/2}p$, but that the subgroup-based maxT may be more powerful if $n^{-1/2}\log p$ is large.
		We confirm our theoretical findings in a simulation study, where we additionally find that the power gap remains large even if the proportion of false hypotheses is moderately small.
		The result is a substantial improvement in power of the group invariance-based MaxT method at only a fraction of the computational cost, compared to the use of the entire group, a large Monte Carlo sample, or a large subgroup.
	
		Finally, we study the $p = 1$ setting under Gaussianity, where the maxT method coincides with testing a single hypothesis.
		Here, we show that the phenomenon cannot occur: the power monotonically increases when passing to a supergroup.
		This explains why the non-monotonicity was not observed by \citet{koning2023more}, since their analysis only considers the $p = 1$ setting.
			
\section{The group invariance-based maxT method}	
	We observe an $n \times p$ matrix 
	\begin{align*}
		X = n^{1/2}\iota\mu' + E,
	\end{align*}
	where $\mu$ is a $p$-vector of means, $\iota$ is a unit vector, and $E$ is a zero-mean error matrix that is invariant under some compact group of orthonormal matrices $\mathcal{G}$.
	That is: $E \overset{d}{=} GE$, for all $G \in \mathcal{G}$.
	
	This includes the group of all permutation matrices, under which invariance is typically called \emph{exchangeability}, but also the \emph{orthogonal group} that consists of all orthonormal matrices, where invariance is often called \emph{sphericity}.
	Another example is the group of sign-flipping matrices, which are diagonal matrices with diagonal elements in $\{-1, 1\}$.
	This group is particularly easy to study as it is finite and commutative, and therefore all its subgroups as well.
	Invariance under the sign-flipping group is equivalent to having symmetric marginal distributions about zero.
	
	The goal is to test the $p$ hypotheses $H_0^j : \mu_j = 0$ against the $p$ alternatives $H_1^j : \mu_j > 0$, $j = 1, \dots, p$, with a method that has high \emph{power}: the expected proportion of correctly rejected hypotheses.
	At the same time, we would like to control the familywise error rate, the probability that at least one of these hypotheses is falsely rejected, by $\alpha \in [0, 1]$.
	As a test statistic, we consider $T : X \mapsto \iota'X$, which can be interpreted as containing $n^{1/2}$ times sample mean of each column of $X$ if $\iota = n^{-1/2}(1, 1, \dots, 1)'$.

	The group invariance-based maxT method works by rejecting $H_0^j$ in favour of $H_1^j$, if $\alpha$ exceeds the $j$th $p$-value $\mathbb{P}_{\overline{G}}\left[\iota'X_j > \max_i \iota'\overline{G}X_i\right]$, where $X_j$ is the $j$th column of $X$, and $\overline{G}$ is uniformly distributed on the group $\mathcal{G}$.
	
	As the group $\mathcal{G}$ is often large, the traditional approach is to instead use a Monte Carlo method based on a random subset.
	In particular, let $\overline{G}_M$ be uniform on $(I, \overline{G}^1, \overline{G}^2, \dots, \overline{G}^{M-1})$, where $I$ is the identity element of $\mathcal{G}$ and the $\overline{G}^i$'s are independent and uniform on $\mathcal{G}$, for all $i$ \citep{hemerik2018exact}.
	The Monte Carlo method then rejects the $j$th hypothesis if $\mathbb{P}_{\overline{G}_M}\left[\iota'X_j > \max_i \iota'\overline{G}_MX_i\right] \leq \alpha$, and still controls the familywise error rate.
	Although we are not aware of a formal proof for the maxT method, one sees in practice that the power is increasing in $M$, with the power of the method based on the entire group as a limit: see also the solid line in the right panel of Figure \ref{fig:intro}.
	For this reason, $M$ is often chosen in the order of hundreds or thousands for common values of $\alpha$, which also ensures good replicability of the method.
	
	\subsection{Selecting a subgroup}
		\citet{dobriban2021consistency} and \citet{koning2023more}  observe that both the Monte Carlo and full-group tests are affected by a ``leak'' from signal into noise, in the sense that the dispersion of the reference distribution increases with the signal $\mu$.
		In the following result, we show that a similar leak appears in the maxT method.
		In particular, the term $n^{1/2}\iota'\overline{G}\iota\mu_1$ distorts the right-hand-side compared to the situation under the null where $\mu = 0$.
	
		\begin{thm}\label{thm:consistency}
			Let $\mathcal{G}$ be a compact group of orthonormal matrices.
			Suppose $E$ is $\mathcal{G}$ invariant with i.i.d. columns, and let $E^2$ be an independent copy of $E$.
			Assume $\mu_1 = \mu_l$, for all $1 \leq l \leq p$.
			Let $M \geq 1$ and let $\alpha = 1/M$.
			Assuming that $p \to \infty$, the Monte Carlo maxT method based on $M$ draws from $\mathcal{G}$ is consistent if and only if 
			\begin{align*}
				\mathbb{P}_{\overline{G}, X_1, E^2}\left[\iota'X_1 > n^{1/2}\iota'\overline{G}\iota\mu_1 +  \max_j \iota'E_j^2\right] \to 1,
			\end{align*}
			where $\overline{G}$ is uniform on $\mathcal{G}$.
		\end{thm}
		
		In the context of testing a single hypothesis, \citet{koning2023more} set out to ``plug'' this leak.
		They note that one still obtains a valid test if $\mathcal{G}$ is replaced by one of its (compact) subgroups, which is a subset that still has a group structure.
		They suggest to carefully select a subgroup with approximately $M$ elements, in order to obtain more power than a Monte Carlo method based on $M$ random samples from the group.
			
		In particular, they argue that subgroups $\mathcal{S}$ for which $\max_{S \in \mathcal{S}\setminus\{I\}}\iota'S\iota$ is minimized are expected to have good power properties.
		They especially highlight ``oracle'' subgroups $\mathcal{S}$ for which the leak vanishes: $\iota'S\iota = 0$, for all $S \in \mathcal{S} \setminus \{I\}$.
		They also note the existence of another type of ``non-positive'' subgroup for which the leak is non-positive: $\iota'S\iota \leq 0$ for all $S \in \mathcal{S} \setminus \{I\}$.
		
		Unfortunately, such oracle and non-positive subgroups only exist up to a size of $n$ and $2n$ respectively, so \citet{koning2023more} recommend to use larger subgroups for which $\max_{S \in \mathcal{S}\setminus\{I\}}\iota'S\iota$ is positive but ``small'', to increase the power that may be lost by using a small subgroup.
		As observed in the left panel of Figure \ref{fig:intro}, and as we prove under Gaussianity in Section \ref{sec:monotonicity}, this recommendation is correct when testing a single hypothesis.
		
		Theorem \ref{thm:consistency} shows that a similar leak appears in the maxT method, so that the same subgroups designed for single hypothesis testing are suitable here as well.
		However, as seen in the right plot on Figure \ref{fig:intro}, the strategy suggested by \citet{koning2023online} to choose a larger subgroup can be detrimental for the power of the maxT method.
		Hence, in contrast to \citet{koning2023more}, we suggest to use oracle and non-positive subgroups for the maxT method, and not larger subgroups.
		These oracle and non-positive subgroups are exactly the subgroups that yield the highest power in Figure \ref{fig:intro}.
		We explore this phenomenon in the following sections.

	\subsection{Monotonicity for a single hypothesis under Gaussianity}\label{sec:monotonicity}
		Theorem \ref{thm:monotonicity} shows that in the special case that $p = 1$ and a Gaussian alternative, a supergroup always yields a more powerful test.
		As a consequence, one should always aspire to use a test based on the entire group $\mathcal{G}$ in this setting.
		\begin{thm}\label{thm:monotonicity}
			Suppose we have a nested sequence of subgroups $\mathcal{S}_1 \subseteq \mathcal{S}_2 \subseteq \dots$.
			Suppose that $p = 1$, and we want to test the null hypothesis that $X$ is $\mathcal{S}_k$ invariant against the alternative hypothesis that $X \sim \mathcal{N}(\iota \mu, \sigma^2 I_n)$, $\mu, \sigma > 0$.
			Then, the $\mathcal{S}_j$ invariance test with test statistic $\iota'X$, is at least as powerful as the $\mathcal{S}_i$ invariance test for $i \leq j \leq k$.
			That is, the power is monotonically increasing in the number of elements of the subgroup.
		\end{thm}
		
		The proof strategy is to show that in this special setting, the $\mathcal{S}_i$ invariance test is a likelihood ratio test for $\mathcal{S}_i$ invariance against Gaussian location shift.
		This follows from an observation in Section 6.2 in \citet{koning2023online}, which generalizes an observation in the final paragraph \citet{lehmann1949theory} for the special case that $\mathcal{S}_i$ is the orthogonal group.
		Applying the Neyman-Pearson lemma, and using the fact that a subgroup invariance test controls size yields the monotonicity property.
		
		As the setting in Theorem \ref{thm:monotonicity} is the primary setting studied by \citet{koning2023more}, the result explains why they did not find any examples in which a subgroup yields a more powerful test. 
		Moreover, the result proves that \citet{koning2023more} were indeed correct in suggesting to use a larger subgroup in this setting.
		However, Theorem \ref{thm:monotonicity} crucially relies on the Neyman-Pearson lemma, so as soon as we move beyond likelihood ratio tests there is little reason to believe that a the monotonicity property would still hold.

\section{Power of the MaxT method under Gaussianity}
	In order to understand the power gap between the subgroup and full-group maxT method, we study their power in a Gaussian location model.

	For oracle subgroups, a particularly clean result about the power of the maxT method can be obtained under Gaussianity, which generalizes a result by \citet{koning2023more} from a single test to the maxT method.
	The result states that if the true distribution is Gaussian, unbeknownst to the analyst, then using an oracle subgroup $\mathcal{S}$ is equivalent to having access to $|\mathcal{S}|$ samples from the true distribution.
	
	\begin{thm}\label{thm:oracle}
		Suppose $E$ has i.i.d. rows from $\mathcal{N}(0, \Sigma)$.
		Suppose that $\mathcal{S}$ is an oracle subgroup.
		Let $\mu$ have $k$ positive elements, and let the remaining elements be equal to zero.
		Let $Z \sim \mathcal{N}(n^{1/2}\mu, \Sigma)$, and $Y$ have i.i.d. rows from $\mathcal{N}(0, \Sigma)$.
		The power of the $\mathcal{S}$-based maxT method is
		\begin{align*}
			\frac{1}{k}\sum_{j=1}^k\mathbb{P}_Z\left[Z_j > q_\alpha^{|\mathcal{S}|}\left(\max_{1\leq l\leq p} Y_l\right)\right],
		\end{align*}
		where $q_\alpha^{|\mathcal{S}|}(\max_{1\leq l\leq p} Y_l)$ denotes the $\alpha$ sample upper-quantile based on $|\mathcal{S}|-1$ draws of $\max_{1\leq l\leq p} Y_l$.
	\end{thm}

	For an arbitrary group $\mathcal{G}$, the power of the maxT method is
	\begin{align*}
		\frac{1}{k} \sum_{j = 1}^k\mathbb{P}_{E}\left[n^{1/2}\mu_j + \iota'E_j > q_\alpha^{\overline{G}}\left(\max_{1 \leq l \leq p} \iota'\overline{G}\iota\mu_l + \iota'\overline{G}E_l\right)\right],
	\end{align*}
	where $\overline{G}$ is uniform on $\mathcal{G}$, and $q_\alpha^{\overline{G}}\left(\max_{1 \leq l \leq p} \iota'\overline{G}\iota\mu_l + \iota'\overline{G}E_l\right)$ denotes the $\alpha$ upper-quantile of the distribution of $\max_{1 \leq l \leq p} \iota'\overline{G}\iota\mu_l + \iota'\overline{G}E_l$ where only $\overline{G}$ is considered random.
	
	Unfortunately, obtaining an insightful characterization of the power as in Theorem \ref{thm:oracle} is \textit{substantially} more involved for non-oracle groups.
	For the special case of the group $\mathcal{H}$ of all orthonormal matrices, independence of the columns of $E$ and $\mu_l = \mu_1$ for all $1 \leq l \leq p$, we empirically find that the power of the maxT method is well approximated by
	\begin{align}\label{eq:full}
		\mathbb{P}_{Z_1}\left[Z_1 > q_\alpha^{\overline{H}, Y}\left(n^{1/2}\iota'\overline{H}\iota\mu_1 + \max_{1 \leq l \leq p} Y_l\right)\right].
	\end{align}
	In Section \ref{appn:cnj} we also provide extensive heuristic arguments for this approximation.
	While these heuristic arguments are of asymptotic nature, simulations suggest the approximation is very accurate for $n, p \geq 10$.
	
	In order to compare the oracle subgroup and full-group approach, we specialize Theorem \ref{thm:oracle} to the setting in \eqref{eq:full} in the following corollary.
	
	\begin{cor}\label{cor:oracle_specialized}
		Suppose $E$ has i.i.d. rows from $\mathcal{N}(0, I)$ and $\mu_l = \mu_1$ for all $1 \leq l \leq k$ and $\mu_l = 0$ for $k < l \leq p$.
		Suppose that $\mathcal{S}$ is an oracle group.
		Let $Z \sim \mathcal{N}(n^{1/2}\mu, I)$ and let $Y$ have i.i.d. rows from $\mathcal{N}(0, I)$.
		Then, the power of the $\mathcal{S}$-based maxT method is
		\begin{align*}
			\mathbb{P}_{Z_1}\left[Z_1 > q_\alpha^{|\mathcal{S}|}\left(\max_{1 \leq l \leq p} Y_l\right)\right].
		\end{align*}
	\end{cor}
		
	Comparing Corollary \ref{cor:oracle_specialized} to \eqref{eq:full}, we notice that for large $p$, $n$ and $|\mathcal{S}|$, the main difference is the appearance of the term $n^{1/2}\iota'\overline{H}\iota\mu_1$ in the reference distribution used for the full-group method, which is approximately $\mathcal{N}(0, \mu_1^2)$-distributed.
	As a consequence, the full-group method's reference distribution is a mean-preserving spread relative to the oracle subgroup method's reference distribution.
	While this is not sufficient to guarantee that all upper-quantiles of the full-group reference distribution are larger, it is typically observed in practice.\footnote{In case the involved distributions are symmetric, a mean-preserving spread is equivalent to the statement that the upper quantiles are larger (Shaked \& Shanthikumar (2007), Section 3.B.1, p. 151.). The distribution of $\iota'\overline{H}\iota n^{1/2}\mu_1$ is symmetric, and while the distribution of $\max_j Y_j$ is not symmetric, it is also not highly asymmetric: it can be approximated by a Gumbel distribution, which has a skewness of approximately 1.14.}
	Furthermore, the mean-preserving spread increases with $\mu_1$, suggesting that the difference between the critical values is expected to be larger if $\mu_1$ is larger.

\section{Relative efficiency under Gaussianity}\label{sec:rel_eff}
		While the previous section provides some insights into why the full-group method can be less powerful than the oracle subgroup-based approach, it does not guarantee that this difference is large in practice.
		Intuitively, while we expect the difference between the critical values to increase in $\mu_1$, the testing problem simultaneously becomes easier as $\mu_1$ increases.
		As a result, it might be that the difference between the critical values only becomes substantial when the power is approximately $1$, which could render the power difference practically insignificant.
		In this section, we show that the power difference \emph{is} substantial in practice, by studying the relative efficiency of the subgroup and full-group-based maxT methods under Gaussianity.
		
		In order to do so, we derive the signals $\mu^{\text{OS}}$ and $\mu^{\mathcal{H}}$ for which the oracle subgroup and full group maxT method have power approximately $1/2$, respectively.
		By studying when the difference between these signals is large, we identify the values of $n$ and $p$ for which the difference between the power of the tests is large.
		The derivations can be found in Section \ref{suppl:derivations} of the Appendix, and the simulations in Section \ref{sec:simulations} demonstrate that the resulting power indeed approaches $1/2$ for sufficiently large $p$ and an appropriately large oracle subgroup.
		
		We find that for sufficiently large $p$ and a sufficiently large subgroup, the oracle subgroup method has power approximately $1/2$ if
		\begin{align*}
			\mu^{\text{OS}}
				&= -n^{-1/2}[\Gamma^{-1}(1-\alpha)/\Phi^{-1}(1/p) + \Phi^{-1}(1/p)],
		\end{align*}
		where $\Phi^{-1}$ and $\Gamma^{-1}$ denote the quantile functions of the standard Gaussian and Gumbel distributions, respectively.
		Moreover, we find that the maxT method based the entire orthogonal group $\mathcal{H}$ has power approximately 1/2 if
		\begin{align*}
			\mu^{\mathcal{H}}
				&= \left(\frac{c^2[a^2 + b(n - c^2)]}{(c^2 - n)^2}\right)^{1/2} + \frac{an^{1/2}}{n - c^2},
		\end{align*}
		where $a = -\gamma / \Phi^{-1}(1/p) - \Phi^{-1}(1/p)$, $b = \pi^2/(6(\Phi^{-1}(1/p))^2)$ and $c = \Phi^{-1}(1-\alpha)$, where $\gamma \approx 0.58$ is the Euler-Mascheroni constant, and $\pi \approx 3.14$ the half-circle constant.
		Using the approximation $\Phi^{-1}(1/p) \approx (2\log p)^{1/2}$, both $\mu^{\text{OS}}$ and $\mu^{\mathcal{H}}$ are order $n^{-1/2}\log^{1/2}p$, so that we expect both methods to be consistent for $\mu \succ n^{-1/2}\log^{1/2}p$.
		
		When comparing $\mu^{\text{OS}}$ and $\mu^{\mathcal{H}}$, we find that to achieve a power of 1/2, the oracle subgroup method requires a smaller signal than the full-group method if
		\begin{align*}
			n^{-1/2} \Gamma^{-1}(1-\alpha) + n^{-1/2}\Phi^{-1}(1/p)^2
				\geq \left[\left(\frac{\gamma - \Gamma^{-1}(1-\alpha)}{\Phi^{-1}(1-\alpha)}\right)^2 - \frac{\pi^2}{6}\right]^{1/2}.
		\end{align*}
		The right-hand side of the inequality does not depend on $p$ or $n$ and is within the range $[0.25, 1.15]$ for $0.01 \leq \alpha \leq 0.1$.
		The left-hand-side is of order $n^{-1/2}\log p$, so that we expect the difference between the tests to be largest in very high-dimensional settings.
		These predictions align with our findings from the simulation study presented in Section \ref{sec:simulations}.
	
\section{Simulation study}\label{sec:simulations}
	In this section, we describe the setup and results of our simulation study.
	For simplicity, we consider the same one-sided Gaussian location setup as studied in the previous sections, with $\iota = (1, \dots, 1)'$, where we aim to maximize the power while controlling the familywise error rate by $\alpha = .05$.
	
	For the subgroup method, we use a sign-flipping oracle subgroup from the R package \url{https://github.com/nickwkoning/NOSdata}.
	For the full-group method, we use the group of all sign-flipping matrices, which are diagonal matrices with $-1$ and $1$ as diagonal elements.
	Since the Gaussian distribution is invariant under the group of sign-flips, the maxT method will control the familywise error rate.
	
	Because the full sign-flipping group is of order $2^n$ and therefore far too large to use in its entirety, we use the Monte Carlo maxT method based on 1000 draws from these groups.
	We consider various values of $p$, $n$ and $\mu$, and also vary the proportion of false hypotheses.
	We also considered a Monte Carlo approach based on the orthogonal group $\mathcal{H}$ that contains all orthonormal matrices, but the resulting power was visually indistinguishable from the method that uses sign-flipping groups.
	
\begin{figure}
		\begin{adjustbox}{width=\textwidth}
		\input{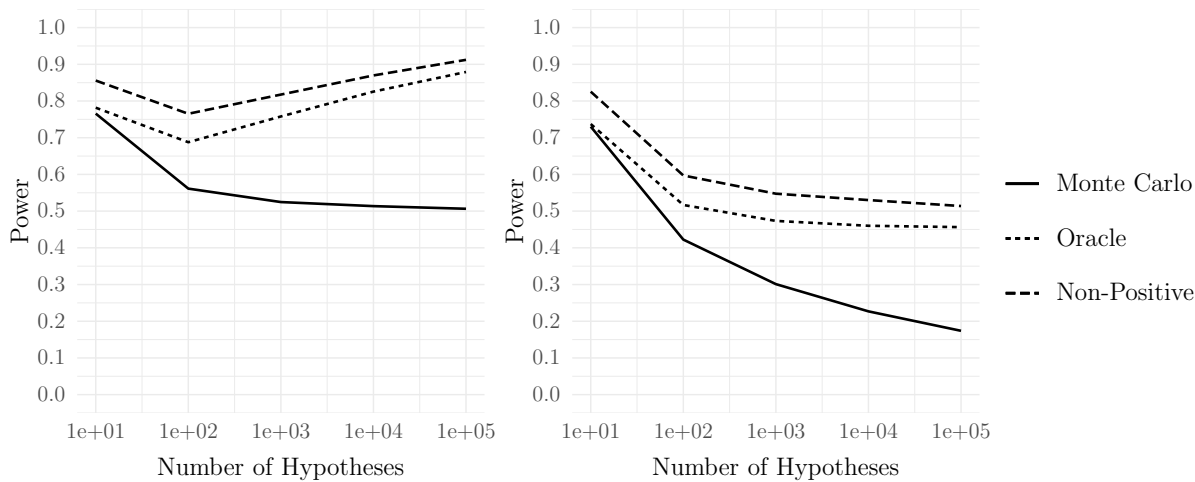}
		\end{adjustbox}
		\caption{The power for differing numbers of hypotheses. The plots are based on 1000 repetitions for $n = 32$, of the Monte Carlo method with 1000 draws (solid), oracle subgroup method with $32$ elements (dotted), non-positive subgroup method with $64$ elements (dashed), using signal $\mu^{\mathcal{H}}$ (left) and $\mu^{\text{OS}}$ (right).}
		\label{fig:muOSvsmuH}
	\end{figure}

	In Figure \ref{fig:muOSvsmuH}, we confirm that the signals $\mu^{\text{OS}}$ and $\mu^{\mathcal{H}}$ from Section \ref{sec:rel_eff} yield a power close to $1/2$ for the oracle subgroup and full-group-based methods, respectively, provided that $p$ is sufficiently large.
	The oracle subgroup seems to yield power slightly below $1/2$ in the right panel of Figure \ref{fig:muOSvsmuH}, which we expect is due to its small size as it consists of just 32 elements.
	Indeed, the non-positive subgroup, which has similar properties but is twice as large, has power very close to $1/2$ for large $p$.
	Moreover, the figure also shows that the power gap is increasing in $p$, as predicted in Section \ref{sec:rel_eff}.
	
\begin{figure}
		\hspace{-1cm}
		\input{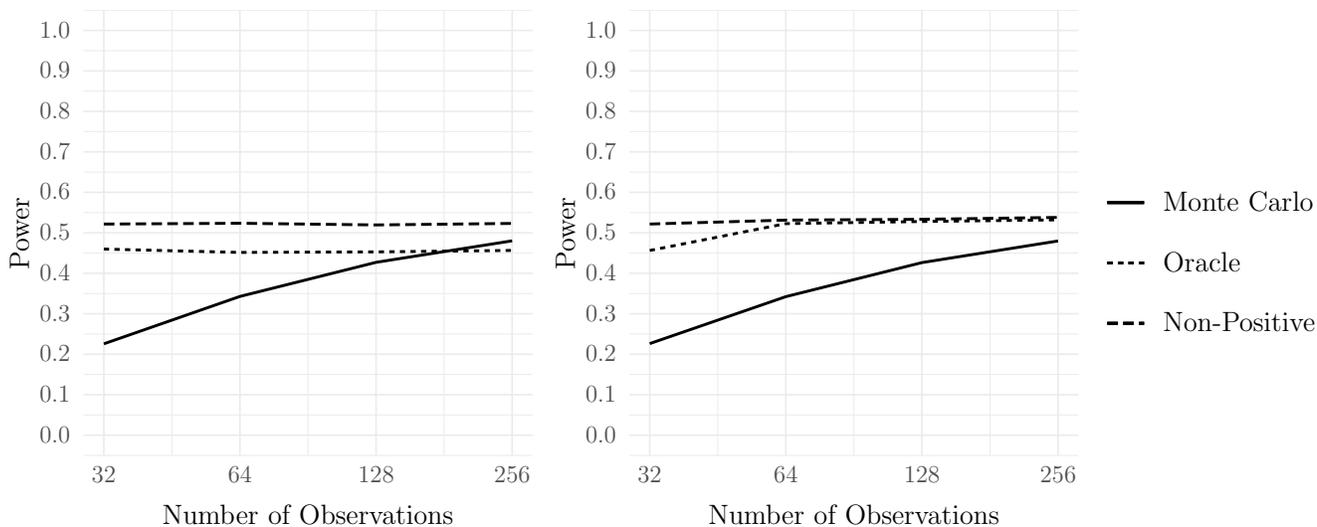}
		\caption{The power for differing numbers of observations. The plots are based on 1000 repetitions for $p = 10000$, of the Monte Carlo method with 1000 draws (solid), oracle subgroup method with $32$ elements (left, dotted) and $n$ elements (right, dotted), non-positive subgroup method with $64$ elements (left, dashed) and $2n$ elements (right, dashed), using signal $\mu^{\text{OS}}$.}
		\label{fig:n}
	\end{figure}

	In Figure \ref{fig:n}, we compare the power of the methods for a signal $\mu^{\text{OS}}$ and varying values of $n$.
	In the left panel, we fix the size of the oracle subgroup at 32 and non-positive subgroup at 64, and in the right panel we set their respective sizes equal to $n$ and $2n$.
	As predicted in Section \ref{sec:rel_eff}, we see that the power gap decreases as $n$ increases.
	Moreover, by comparing both panels we observe that the gap between the oracle and non-positive subgroup methods seems independent of $n$, but closes as the size of the subgroup increases.
	This suggests the power gap between the two subgroup methods is mainly due to the size of the subgroup.
	
\begin{figure}
		\hspace{-1cm}
		\input{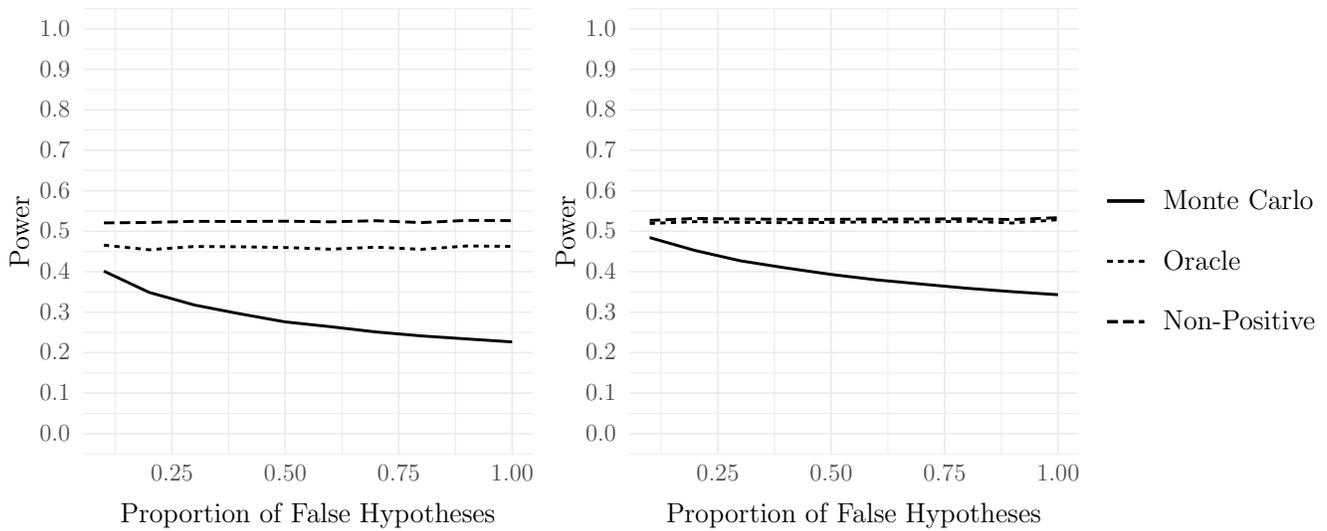}
		\caption{The power for differing proportions of false hypotheses. The plots are based on 1000 repetitions for $p = 10000$, $n = 32$ (left) and $n = 64$ (right), of the Monte Carlo method with 1000 draws (solid), oracle subgroup method with $n$ elements (dotted), non-positive subgroup method with $2n$ elements (dashed), using signal $\mu^{\text{OS}}$.}
		\label{fig:false_hyp}
	\end{figure}

	In Figure \ref{fig:false_hyp}, we vary the proportion of false hypotheses.
	In the left panel, we use $n = 32$ and in the right panel we use $n = 64$.
	We observe that the power gap increases in the proportion of false hypotheses, but that the power difference remains substantial even if the proportion of false hypotheses is moderately small.
	Moreover, the subgroup-based tests seem unaffected by the proportion of false hypotheses, as predicted by Theorem \ref{thm:oracle}.		

\section{Discussion}		
	We believe a fruitful direction for future research the improvement of other permutation- or group-based methods by the strategic selection of a subgroup.
	Moreover, \citet{Ramdas2023} recently suggested permutation tests based on arbitrary distributions, which goes beyond (uniform distributions on) subgroups.
	With this additional flexibility, it may be possible to extract even more power, though it remains unclear how such a distribution can be selected.
		
\section{Acknowledgements}
\if1\blind
{
	We thank Jesse Hemerik and Stan Koobs for useful comments.
} \fi

\if0\blind
{

} \fi

	We have no funding nor any conflicts of interest to disclose.

\clearpage
	
\clearpage

\begin{center}
	{\Large\bf Appendix}
\end{center}

\section{Proof of Theorem \ref{thm:consistency}}\label{appn:consistency}
	\begin{proof}
		As the elements of $\mu$ are equal, we have that the power equals
		\begin{align*}
			\frac{1}{p}\sum_{j=1}^p\mathbb{P}_{\mathcal{S}_M, E}\left[n^{1/2}\mu_1 + \iota' E_j > q_\alpha^{\overline{S}_M}\left(n^{1/2}\iota'\overline{S}_M\iota\mu_l + \max_{l \in \{1, \dots, p\}} \iota'\overline{S}_ME_l\right)\right].
		\end{align*}
		As the columns of $E$ are i.i.d., the probability that $E_j$ is the maximizer equals $1/p$.
		Hence, the difference between 
		\begin{align*}
			\frac{1}{p}\sum_{j=1}^p\mathbb{P}_{\mathcal{S}_M, E}\left[n^{1/2}\mu_1 + \iota' E_j > q_\alpha^{\overline{S}_M}\left(\iota'\overline{S}_M\iota\mu_1 + \max_{l \in \{1, \dots, p\}}\iota'\overline{S}_ME_l\right)\right]
		\end{align*}
		and
		\begin{align*}
			\frac{1}{p}\sum_{j=1}^p\mathbb{P}_{\mathcal{S}_M, E}\left[n^{1/2}\mu_1 + \iota' E_j > q_\alpha^{\overline{S}_M}\left(\iota'\overline{S}_M\iota\mu_1 + \max_{l \in \{1, \dots, p\}\setminus\{j\}}\iota'\overline{S}_ME_l\right)\right]
		\end{align*}
		vanishes as $p \to \infty$.
		Moreover, as the columns of $E$ are i.i.d and choosing $j = 1$ without loss of generality, the latter equals
		\begin{align*}
			\mathbb{P}_{\mathcal{S}_M, E}\left[n^{1/2}\mu_1 + \iota' E_1 > q_\alpha^{\overline{S}_M}\left(\iota'\overline{S}_M\iota\mu_1 + \max_{l \in \{2, \dots, p\}}\iota'\overline{S}_ME_l\right)\right],
		\end{align*}
		which for large $p$ is close to 
		\begin{align*}
			\mathbb{P}_{\mathcal{S}_M, E_1, E^2}\left[n^{1/2}\mu_1 + \iota' E_1 > q_\alpha^{\overline{S}_M}\left(\iota'\overline{S}_M\iota\mu_1 + \max_{l \in \{1, \dots, p\}}\iota'\overline{S}_ME_l^2\right)\right],
		\end{align*}
		where $E^2$ is an independent copy of $E$.
		Hence, assuming that $p \to \infty$, the maxT method is consistent if and only if the term in the previous display converges to $1$.

		Let $\overline{S}_M$ be uniform on $\mathcal{S}_M = (I, \overline{S}^1, \overline{S}^2, \dots, \overline{S}^M)$, where $I$ is the identity element of $\mathcal{S}$ and the $\overline{S}^i$s are independent and uniform on $\mathcal{S}$, for all $i$.
		Let $\overline{S}$ independently be uniform on $\mathcal{S}$.
		Analogous to Lemma 4.2 in the Supplementary Material of \citet{dobriban2021consistency}, we have that 
		\begin{align*}
			&\mathbb{P}_{\mathcal{S}_M, E_1, E^2}\left[n^{1/2}\mu_1 + \iota' E_1 > q_\alpha^{\overline{S}_M}\left(\max_j n^{1/2}\iota'\overline{S}_M\iota\mu_j + \iota'\overline{S}_ME_j^2\right)\right] \\
				&= \mathbb{P}_{\mathcal{S}_M, E_1, E^2}\left[n^{1/2}\mu_1 + \iota' E_1 > q_{1/|\mathcal{S}|}^{\overline{S}_M}\left(\max_j n^{1/2}\iota'\overline{S}_M\iota\mu_j + \iota'\overline{S}_ME_j^2\right)\right] \\
				&= \mathbb{P}_{\mathcal{S}_M, E_1, E^2}\left[n^{1/2}\mu_1 + \iota' E_1 > \max_{S \in \mathcal{S}_M}\max_j n^{1/2}\iota'S\iota\mu_j + \iota'SE_j^2\right] \\
				&= \mathbb{E}_{E_1, E^2} \mathbb{P}_{\overline{S}}\left[n^{1/2}\mu_1 + \iota' E_1 > \max_j n^{1/2}\iota'\overline{S}\iota\mu_j + \iota'\overline{S}E_j^2\right]^{M},
		\end{align*}
		which converges to 1 if and only if
		\begin{align*}
			\mathbb{P}_{\overline{S}, E_1, E^2}\left[n^{1/2}\mu_1 + \iota' E_1 > \max_j \iota'\overline{S}\iota\mu_j + \iota'\overline{S}E_j^2\right] 
				\to 1.
		\end{align*}
		As $E^2$ is $\mathcal{G}$ invariant, we have $E^2 \overset{d}{=} GE^2$ for all $G \in \mathcal{G}$.
		Then, as $\mathcal{S} \subseteq \mathcal{G}$, we have $\overline{S} \in \mathcal{G}$.
		Hence, $\overline{S}E^2 \overset{d}{=} E^2$. 
		As a result,
		\begin{align*}
			\hspace{-.7cm}
			\mathbb{P}_{\overline{S}, E_1, E^2}\left[n^{1/2}\mu_1 + \iota' E_1 > \max_j n^{1/2}\iota'\overline{S}\iota\mu_j + \iota'\overline{S}E_j^2\right] 
				&= \mathbb{P}_{\overline{S}, E_1, E^2}\left[n^{1/2}\mu_1 + \iota' E_1> \max_j n^{1/2}\iota'\overline{S}\iota\mu_j + \iota'E_j^2\right],
		\end{align*}
		which finishes the proof.
	\end{proof}

	\begin{rmk}
		While we assume in Theorem \ref{thm:consistency} that the columns of $E$ are independent and the elements of $\mu$ are equal, this is mostly out of mathematical convenience.
		We expect analogous results to hold for specific dependence structures of $E$ as long as the columns of $E$ are not too dependent, and the elements of $\mu$ are not too different.
	\end{rmk}

\section{Proof of Theorem \ref{thm:monotonicity}}\label{appn:monotonicity}
	\begin{proof}
		By an observation in Section 6.2 in \citet{koning2023online}, we have that the $\mathcal{S}_i$ invariance test is a likelihood ratio test for testing $\mathcal{S}_i$ invariance against $\mathcal{N}(\iota\mu, \sigma^2I_n)$, for every $i$.
		By the Neyman-Pearson lemma, this test is uniformly most powerful.
		As invariance under a group implies invariance under each of its subgroups, a $\mathcal{S}_i$ invariance test also controls size for $i \leq k$.
		As a result, we can immediately conclude that the $\mathcal{S}_k$ invariance test is at least as powerful a the $\mathcal{S}_i$ invariance test for $i \leq k$.
		
		To show the monotonicity, it remains to show that a $\mathcal{S}_j$ invariance test is at least as powerful as a $\mathcal{S}_i$ invariance test for testing $\mathcal{S}_k$ invariance.
		As $\mathcal{S}_k$ invariance implies $\mathcal{S}_j$ invariance and $\mathcal{S}_i$ invariance, $i \leq j \leq k$, we can also consider testing $\mathcal{S}_j$ invariance.
		By the Neyman-Pearson lemma, the $\mathcal{S}_j$ invariance test is most powerful, and so more powerful than the $\mathcal{S}_i$ invariance test.
		This proves the claim.
	\end{proof}

\section{Proof of Theorem \ref{thm:oracle}}
	The proof strategy of Theorem \ref{thm:oracle} mimics to that of Theorem 8 in \citet{koning2023more}, and generalizes their result from a single hypothesis test to the maxT method.
	\begin{proof}
		First, assume that $|\mathcal{S}| = n$.
		Define the matrix $\mathfrak{S} = \{S\iota\ |\ S \in \mathcal{S}\}$, such that its first column is $\iota$.
		Using $e_1 = (1, 0, \dots, 0)'$, we have
		\begin{align*}
			\mathfrak{S}'X 
				= \mathfrak{S}'\iota\mu' + \mathfrak{S}'E
				= e_1\mu' + \mathfrak{S}'E
				\overset{d}{=} e_1\mu' + E,
		\end{align*}
		where the second equality follows from the fact that $\mathfrak{S}$ is orthonormal and has first column $\iota$, and the equality in distribution from the orthogonal invariance of $E$.
		As a consequence, the second to the final rows of $\mathfrak{S}'X$ have distribution $\mathcal{N}(0, \Sigma)$.
		
		The $\mathcal{S}$-based maxT method rejects the $j$th hypothesis if the $j$th element of the first row of the matrix $(e_1\mu' + E)$ is larger than the row-wise maxima of its remaining rows.
		Since all but the first row of $e_1\mu'$ are equal to zero, the remaining rows of $(e_1\mu' + E)$ equal those of $E$.
		The result for the $|\mathcal{S}| = n$ follows from noting that the rows of $E$ are i.i.d. $\mathcal{N}(0, \Sigma)$-distributed.
		
		For the $|\mathcal{S}| \leq n$ case, $\mathfrak{S}$ is $n \times |\mathcal{S}|$.
		Enlarging $\mathfrak{S}$ by adding $(n - |\mathcal{S}|)$ columns of zeros, the above reasoning can be extended to prove the result for oracle subgroups of arbitrary size.
	\end{proof}

\section{Heuristic arguments for equation \eqref{eq:full}}\label{appn:cnj}
	In \eqref{eq:full} we presented the following conjecture.
	
	\begin{cnj}\label{cnj:full}
		Define $\overline{h} = \overline{H}\iota$.
		Suppose $E$ has i.i.d. rows $\mathcal{N}(0, I)$ and $\mu_l = \mu_1$ for all $1 \leq l \leq p$.
		Let $Z_1 \sim \mathcal{N}(n^{1/2}\mu_1, 1)$ and let $Y$ have i.i.d. rows $\mathcal{N}(0, I)$.
		As $p, n \to \infty$, we have
		\begin{align*}
			\hspace{-1cm}
			\left|\frac{1}{p}\sum_{j=1}^p \mathbb{P}_E \left[n^{1/2}\mu_1 + \iota'E_j > q_\alpha^{\overline{h}}\left(n^{1/2}\overline{h}'\iota\mu_1 + \max_i \overline{h}'E_i\right) \right]- \mathbb{P}_{Z_1}\left[Z_1 > q_\alpha^{\overline{h}, Y}\left(n^{1/2}\overline{h}'\iota\mu_1 + \max_{1 \leq l \leq p} Y_l\right)\right]\right|
			\to 0.
		\end{align*}
	\end{cnj}
		
	Our heuristic arguments for why we believe this conjecture holds consists of several steps.
	First, as $p \to \infty$, we can use reasoning analogous to the start of the proof of Theorem \ref{thm:consistency}, the independence of the columns of $E$ and the fact that all elements of $\mu$ are equal to find that for large $p$
	\begin{align*}
		\frac{1}{p}\sum_{j = 1}^{p}\mathbb{P}_{E}\left[n^{1/2}\mu_1 + \iota'E_j > q_\alpha^{\overline{h}}\left(\overline{h}'\iota n^{1/2}\mu_1 + \max_i \overline{h}'E_i\right)\right].
	\end{align*}
	is close to 
	\begin{align*}
		\mathbb{P}_{E, Z}\left[Z > q_\alpha^{\overline{h}}\left(\overline{h}'\iota n^{1/2}\mu_1 + \max_i \overline{h}'E_i\right)\right],
	\end{align*}	
	where $Z \sim \mathcal{N}(n^{1/2}\mu_1, 1)$.
	It remains to argue that the conditional distribution of $\overline{h}'\iota n^{1/2}\mu_1 + \max_i n^{1/2}\overline{h}'E_i$, given $E$, is close to the distribution of $\overline{h}'\iota n^{1/2}\mu_1 + \max_i Y_i$ with high probability.
	
	To do this, we start by arguing that the conditional distribution of $\max_i n^{1/2}\overline{h}'E_i$ given $E$ is close to the distribution of $\max_i Y_i$, with high probability.
	The argument consists of three components:
	\begin{enumerate}
		\item By Theorem 3 in \citet{jiang2006how}, we have that $E$ is well-approximated by $n^{1/2}\Gamma$, where $\Gamma$ is the $n \times p$ matrix with orthonormal columns obtained by applying the Gram-Schmidt orthogonalization to $E$.
		\item By the orthogonal invariance of $\overline{h}$, we have that $n^{1/2}\overline{h}'\Gamma$ is equal in distribution to $n^{1/2}\overline{h}_{[1:p]}'$, where $\overline{h}_{[1:p]}$ contains the first $p$ elements of $\overline{h}$.
		\item By Theorem 1 in \citet{jiang2006how}, $n^{1/2} \overline{h}_{[1:p]} \overset{d}{\to} \mathcal{N}(0, I_p)$. Hence, for large $n$ and $p$, the distribution of $\max_i n^{1/2} \overline{h}_{[1:p]}$ is close to that of $\max_i Y_i$.
	\end{enumerate}	
	
	It remains to show that $\iota'\overline{h}n^{1/2}\mu_1$ is asymptotically independent from $\max_i n^{1/2}\overline{h}E_i$. 
	To do so, we can go back to step 1 and add $\iota$ as a column to $E$ and apply the Gram-Schmidt orthognalization.
	Following steps 2 and 3 then yields an additional independent standard normal element.
	
	Unfortunately, Theorem 3 in \citet{jiang2006how} requires $p = o(n / \log n)$ and Theorem 1 in \citet{jiang2006how} requires $p = o(n^{1/2})$ and $p \leq n$, which both exclude the high-dimensional regime we are interested in.
	As far as we are aware, there currently do not exist tools that can satisfactorily deal with the high-dimensional regime. 
	Intuitively, the above strategy breaks down due to fundamental restriction that an $n$-dimensional vector can have at most $n$ orthogonal elements, and at most $n$ independent marginals.
	
	One potential route to escape this fundamental restriction is an asymptotic independence representation by a phantom distribution (see e.g. \citet{jakubowski1993asymptotic}).
	The key idea is that while $p$ identically distributed but ``weakly'' dependent random variables may asymptotically be poorly approximated by $p$ i.i.d. random variables, the \textit{maximum} of the $p$ weakly dependent random variables may still be well-approximated by the \textit{maximum} of $p$ i.i.d. random variables.
	Unfortunately, as far as we are aware, this tool has only been developed under sequential dependence.
	The dependence structure we face is of a different nature: conditional on $E$, the elements of $\overline{h}'E$ are dependent, but the dependence weakens with high $\mathbb{P}_{E}$-probability as $n\to\infty$ (see \citet{cai2013distributions}).
	While it seems possible to extend the results to this type of dependence structures, we consider this (far) beyond the scope of the current paper.

\section{Derivations relative efficiency}\label{suppl:derivations}
	By Corollary \ref{cor:oracle_specialized} and the symmetry of the normal distribution, we know the oracle subgroup-based maxT method has power $1/2$ if $Z_1 = q_\alpha^{|\mathcal{S}|}\left(\max_j Z_j\right)$.
	For large $p$, $\max_j Z_j \approx \text{Gumbel}(-\Phi^{-1}(1/p), -1/\Phi^{-1}(1/p))$, where $\Phi^{-1}$ is the quantile function of the standard normal distribution.
	Hence, if $|\mathcal{S}|$ is sufficiently large, we expect the oracle subgroup-based maxT method to have power $\approx 1/2$ if
	\begin{align*}
		\mu^{\text{OS}}
			&= -n^{-1/2}[\Gamma^{-1}(1-\alpha)/\Phi^{-1}(1/p) + \Phi^{-1}(1/p)] \\
			&\approx \log (1/\alpha) (2n\log p)^{-1/2} + n^{-1/2} (2\log p)^{1/2}.
	\end{align*}
	Similarly, by \eqref{eq:full}, the Monte Carlo test has power approximately 1/2 if $\mu^{\mathcal{H}}$ solves
	\begin{align}\label{eq:to_solve_full_group}
		n^{1/2}\mu_1^{\mathcal{H}} 
			= q_\alpha\left(\iota'\overline{H}\iota n^{1/2} \mu_1^{\mathcal{H}} + \max_{1 \leq l \leq p} Y_l\right).
	\end{align}
	Besides the Gumbel approximation, we can use that $\sqrt{n}\iota'\overline{H}\iota \approx \mathcal{N}(0, 1)$, even for quite small $n$ (see e.g. chap. 7 in \citet{eaton1989group}), so that
	\begin{align*}
		\iota'\overline{H}\iota n^{1/2} \mu_1^{\mathcal{H}} + \max_{1 \leq l \leq p} Y_l
			\overset{\text{approx}}{\sim} \mathcal{N}(0, \mu_1^{\mathcal{H}}) + \text{Gumbel}(-\Phi^{-1}(1/p), -1/\Phi^{-1}(1/p)),
	\end{align*}
	where the $+$ on the right-hand-side denotes a convolution of the probability measures.
	Unfortunately, we could not find any work on the convolution of a normal and gumbel random variable.
	However, if $\mu_1^{\mathcal{H}}$ is not too small, it seems well approximated by the normal distribution
	\begin{align*}
		\mathcal{N}(-\gamma / \Phi^{-1}(1/p) - \Phi^{-1}(1/p), \pi^2/6(\Phi^{-1}(1/p))^2 + (\mu_1^{\mathcal{G}})^2).
	\end{align*}
	Substituting this approximation into equation \eqref{eq:to_solve_full_group} and solving for $\mu_1^{\mathcal{H}}$ yields
	\begin{align*}
		\mu_1^{\mathcal{G}}
			&= \sqrt{\frac{c^2(a^2 + b(n - c^2))}{(c^2 - n)^2}} + \frac{a\sqrt{n}}{n - c^2},
	\end{align*}
	where $a = -\gamma / \Phi^{-1}(1/p) - \Phi^{-1}(1/p)$, $b = \pi^2/(6(\Phi^{-1}(1/p))^2)$ and $c = \Phi^{-1}(1-\alpha)$.
	
	Substituting $\mu^{\text{OS}}$ into \eqref{eq:to_solve_full_group} and solving for the terms that involve $p$ and $n$, we find that the MC test has power less than 1/2 when the OS method has power 1/2 if
	\begin{align*}
		\left[\left(\frac{\gamma-\Gamma^{-1}(1-\alpha)}{\Phi^{-1}(1-\alpha)}\right)^2 - \frac{\pi^2}{6}\right]^{1/2}
			\geq n^{-1/2}\Phi^{-1}(1/p)^2 + n^{-1/2}\Gamma^{-1}(1-\alpha).	
	\end{align*}
	
	To verify this, we also check when the OS test has power at least $1/2$ when the full-group test has power $1/2$, by comparing $\mu^{\mathcal{H}}$ to $\mu^{\text{OS}}$.
	Using $(c^2 - n)^2 \approx n^2$, we have
	\begin{align*}
		 \mu_1^{\text{OS}}
		 	\geq \mu_1^{\mathcal{H}},
	\end{align*}
	iff
	\begin{align*}
		\left[\left(\frac{\gamma - \Gamma^{-1}(1-\alpha)}{\Phi^{-1}(1-\alpha)}\right)^2 - \frac{\pi^2}{6}\right]^{1/2}
			\geq n^{-1/2}\Phi^{-1}(1/p)^2 + n^{-1/2} \gamma.
	\end{align*}
	Notice that this yields almost exactly the same result, but with $\Gamma^{-1}(1-\alpha)$ replaced by $\gamma$.
	
	\clearpage
	\bibliographystyle{agsm}
	\bibliography{Bibliography-MM-MC}
\end{document}